\documentclass[11pt]{amsart}
\usepackage{amssymb,amsmath}
\usepackage{hyperref}
\usepackage{bm}
\usepackage[]{graphicx}
\usepackage{amssymb, epsfig}
\numberwithin{equation}{section}

\numberwithin{figure}{section}
\newtheorem{theorem}{Theorem}[section]
\newtheorem{lemma}[theorem]{Lemma}

\newtheorem{remark}[theorem]{Remark}
\begin{document}
%
%
%
%
\title[An elementary observation on the Robin problem for NLS]{An elementary observation on the Robin problem for NLS}

\author[D.~C. Antonopoulou]{D.~C. Antonopoulou$^{\$}$}

\author[A. Chatziafratis]{A. Chatziafratis$^{**}$}

\author[S. Kamvissis]{S. Kamvissis$^{*}$}
\thanks{$^*$ Mathematics Building,
University of Crete, GR--700 13 Heraklion, Greece, and, Institute
of Applied and Computational Mathematics, FORTH, GR--711 10
Heraklion, Greece, email: spyros@tem.uoc.gr}
\thanks{$^\$$ Department of Mathematics, University
of Athens, Greece, email: danton@math.uoa.gr}
\thanks{$^{**}$ Waseda Research Institute for Science and Engineering, Waseda University, Tokyo, Japan, email: chatziafrati@math.uoa.gr}

%
%
\begin{abstract}
We present an infinite set of explicit solutions of the integrable focusing NLS equation on the half-line that are rapidly decaying and satisfy a homogeneous Robin condition,
but have very different behaviors on the wall $x=0$. We then note  a similar phenomenon for the non-homogeneous Robin problem. We argue that this renders a potential application
of the unified transform method to non-linearizable Robin problems impossible.
%
\end{abstract}
%
%
\maketitle
\pagestyle{myheadings}
\thispagestyle{plain}
%
%
%
\section{NLS  and the $N$-soliton}

We consider the focusing NLS equation, with given initial data and a simple inhomogeneous Robin condition:
\begin{equation}\label{eq:nls}
iu_t+u_{xx}+2|u|^2u=0\quad(x>0,\ t>0),\qquad u_x(0,t)+q\,u(0,t)=h(t),\qquad u(x,0)=u_0(x),
\end{equation}
with $q\in\mathbb R$.
The homogeneous case $h(t)=0$ is famously known to be integrable and  have soliton solutions. We present  a specific formula below.
\subsection{The full-line $N$-soliton}\label{sec:nsol}
Let $\zeta_j=i\eta_j$ with $\eta_1>\dots>\eta_N>0$, and let $c_j\in\mathbb R\setminus\{0\}$. Put
\[
\theta_j=-i\zeta_jx-2i\zeta_j^2t,\qquad v_j=\big(c_je^{\theta_j},\,e^{-\theta_j}\big),\qquad
M_{jk}=\frac{\overline{v_{j1}}v_{k1}+\overline{v_{j2}}v_{k2}}{\zeta_k-\overline{\zeta_j}} ,
\]
and define the reflectionless $N$-breather \cite{Y}, \cite{Y2}
\begin{equation}\label{eq:nsol}
U(x,t)=2i\sum_{j,k=1}^N v_{j1}\,(M^{-1})_{jk}\,\overline{v_{k2}} .
\end{equation}

For $N=1$ this is $U(x,t)=2 sgn(c) \cdot \eta sech(2\eta x-\log|c|)\,e^{4i\eta^2t}$.

The numbers $\zeta_j$ are of course the eigenvalues of the associated Zakharov-Shabat operator. By choosing them purely imaginary we ensure that the solitons have zero speed and lie "on top of each other".
They do oscillate in time, hence the name "breather".

Since $\zeta_k-\overline{\zeta_j}=i(\eta_j+\eta_k)$, $iM$ is a Gram matrix with kernel $1/(\eta_j+\eta_k)$, which is positive definite. Hence $\det M\neq0$ for every $x$, $t$.

\subsection{The half-line $N$-soliton}
Restricted to $x>0$, $U$ solves \eqref{eq:nls} with $h=0$ in the following cases \cite{BB}
\begin{itemize}
\item $q=0$ and $c_j=\pm 1$; the Robin condition degenerates into a Neumann condition.
\item $q\neq0$, $|q|<2\eta_j$ for all $j$, with
\[
c_j=\pm\Big(\frac{2\eta_j+q}{2\eta_j-q}\Big)^{(-1)^{N+1}/2}.
\]

\end{itemize}

The trace $U(0,\cdot)$ is exactly quasi-periodic, with frequencies in the $\mathbb Z$-module generated by $\omega_1,\dots,\omega_N$, where $\omega_j=2\eta_j^2$.
Indeed, at $x=0$, $\theta_j=2i\eta_j^2t$, so $v_{j1}=c_je^{i\omega_jt}$ and $v_{j2}=e^{-i\omega_jt}$. Every entry of $M$ and every factor in \eqref{eq:nsol} is a trigonometric polynomial in the angles $\phi_j=\omega_jt$. The matrix $iM$ is positive definite for arbitrary unimodular phases, so $\det M\neq0$. Hence $U(0,t)$ is a rational function of the $e^{i\phi_j}$ with non-vanishing denominator.

There is no decaying radiation term \footnote{There are solutions of the Robin problem \cite{LL} where  the Dirichlet trace is only asymptotically quasi-periodic.} because the data are reflectionless. So the quasi-periodicity is exact for all $t\ge0$. If the $\omega_j$ are rationally independent, the trace has $N$ independent frequencies. Since the $\eta_j$ are arbitrary, there are different smooth, exponentially decaying data that give different frequency sets of any size.

\begin{remark}
The above elementary example shows that (in any reasonable sense) the map taking the data  for the IBVP (initial data and the constant $k$) to the Dirichlet trace is very unstable in terms of long time asymptotics.
Small changes in (the parameters of) the initial data can give non-decaying changes in the Dirichlet trace.
\end{remark}

\section{The non-homogeneous Robin problem}\label{sec:forced}

If $h$ is not zero the problem is non-linearizable in the sense of \cite{FIS}. The question still arises whether it is still  in a sense integrable.
That might be so if the unified transform method introduced by Fokas \cite{F} could transform the IBVP to a Riemann-Hilbert problem
with data that depend continuously on the data of the IBVP and that is amenable to a steepest descent asymptotic analysis. This would be impossible
however if the instability mentioned above persists,  small changes in the initial data giving rise to many different non-decaying changes in the Dirichlet trace; simply because
any possible scattering theory for the boundary data requires a fixed asymptotic t-behavior. Here on the contrary we have a continuum of different asymptotically quasi-periodic
behaviors.

This is indeed what the following perturbation analysis shows.

An elementary calculation gives the time derivatives of mass $ \mathcal N=\int_0^\infty|u|^2\,dx,$ and the energy
$ \mathcal E_q=\int_0^\infty\big(|u_x|^2-|u|^4\big)dx-q\,|u(0,t)|^2$ in terms of  the non-homogeneous term $h$.

\begin{lemma}\label{lem:bal}
Let $u$ be a smooth solution of \eqref{eq:nls}, decaying as $x\to\infty$. Then \begin{align}
\frac{d}{dt}\mathcal N&=2\operatorname{Im}\big(\overline{u(0,t)}\,h(t)\big),\\
\frac{d}{dt}\mathcal E_q&=-2\operatorname{Re}\big(\overline{u_t(0,t)}\,h(t)\big).
\end{align}
\end{lemma}

\begin{proof}
\emph{Mass.} From $\partial_t|u|^2=-2\partial_x\operatorname{Im}(\bar uu_x)$ we get $\dot{\mathcal N}=2\operatorname{Im}(\bar uu_x)(0)$. Substituting $u_x(0)=-qu(0)+h$, the $q$-term is real and drops out.

\emph{Energy.} Integrating by parts,
\[
\frac{d}{dt}\int_0^\infty\big(|u_x|^2-|u|^4\big)=-2\operatorname{Re}(\bar u_tu_x)(0)-2\operatorname{Re}\int_0^\infty\bar u_t\big(u_{xx}+2|u|^2u\big).
\]
The last integral is $-2\operatorname{Re}\int\bar u_t(-iu_t)=0$. At the boundary,
\[
-2\operatorname{Re}(\bar u_tu_x)(0)=2q\operatorname{Re}(\bar u_tu)(0)-2\operatorname{Re}(\bar u_th)=q\,\tfrac{d}{dt}|u(0,t)|^2-2\operatorname{Re}(\bar u_th),
\]
which gives (2.2).
\end{proof}

For $h=0$ both quantities are of course conserved.

\subsection{Values on forced solitons}

When a non-homogenous $h$ term appears, we can still consider the evolution of the eigenvalues of the Zakharov-Shabat operator.
They of course can move in time.

\begin{theorem}\label{prop:val}
For the Robin-compatible $1$-soliton of the first section with $A=2\eta$,
\begin{equation}\label{eq:values}
\mathcal N=A-q,\qquad \mathcal E_q=-\tfrac13\Big(A^3-q^3\Big).
\end{equation}
\end{theorem}

\begin{proof}
Write $u=A sech(Ax+c)e^{iA^2t}$ with $s=tanh c=q/A$. Then
\[
\mathcal N=A(1-s),\qquad \int_0^\infty|u_x|^2=\tfrac{A^3}3(1-s^3),\qquad \int_0^\infty|u|^4=A^3\big(\tfrac23-s+\tfrac{s^3}3\big),\qquad |u(0)|^2=A^2(1-s^2).
\]
Combining these gives $\mathcal N=A-q$ and $\mathcal E_q=-\tfrac13(A^3-q^3)$.
\end{proof}

The above formulae easily generalize to the case of several solitons.

\subsection{First-order frequency shifts}
Let $h=\varepsilon h_1$ with $h_1$ integrable and rapidly decaying, and let the unforced solution be a "wall" $1$-soliton $U$. Assume for the moment
that the forced solution will be close to a  $1$-soliton with amplitude $A+\delta A$, plus radiation. \footnote{The radiation carries only $O(\varepsilon^2)$ of mass and energy.}

Integrate (2.1) and (2.2) with the unforced trace inserted, and use \eqref{eq:values}. To first order,
\begin{equation}\label{eq:first}
\delta A=2\operatorname{Im}\int_0^\infty\overline{U(0,t)}\,h(t)\,dt,\qquad
A^2\,\delta A=2\operatorname{Re}\int_0^\infty\overline{U_t(0,t)}\,h(t)\,dt,
\end{equation}
and $\delta\omega=2A\delta A$.

Since $U(0,t)=\sqrt{A^2-q^2}\,e^{iA^2t}$ and $U_t=iA^2U$,  the two relations in \eqref{eq:first} coincide. They give
\[
\delta A=2\sqrt{A^2-q^2}\,\operatorname{Im}\int_0^\infty h(t)e^{-iA^2t}\,dt .
\]

This calculation can be generalized to the case of $N$ solitons, where we will end up for a system for the different  $\delta\omega_j$ and $\delta A_j$. In the general case the energy
calculation becomes relevant.  All  $\delta\omega_j$ and $\delta A_j$ are real.

One reason to assume that a small perturbation does not destroy solitons is the following. All solitons corresponding to the eigenvalues $i\eta_j$  correspond to zeros
of some scattering data that depend $continuously$ on the perturbation. These could move somewhat under a small enough perturbation, but cannot disappear. For this to happen they would have to hit the real line;
but all zeros are at least at distance $q$ away from the real line. 

Now, because of symmetries ($\zeta \to -\bar \zeta$), the eigenvalues $have~to$ also remain on the imaginary line, unless there is a collision of two imaginary eigenvalues. This does not happen for
a small enough perturbation. Breathers remain breathers even if frequencies can change.  

Note also that even if when there is collision of two imaginary eigenvalues that gives rise to eigenvalues with non-zero real parts, the leading order asymptotics at $x=0$ still come from the remaining
eigenvalues with zero real parts!

If $h$ has compact support and the perturbation size is small, clearly one can ensure that eigenvalues remain imaginary. For decaying $h$ not of compact support, some
eigenvalues can eventually leave the imaginary line. But in general not all of them will.

To summarize, we end up with a situation similar to the unforced case: continua of solutions with different non-decaying asymptotically quasi-periodic
asymptotics at the wall $x=0$. This renders a unified transform impossible!

\subsection{Numerical tests}

We have asked Claude (Anthropic) to provide numerical runs where indeed small non-homogeneities just shift amplitudes and frequencies of breather solutions. 
We mention three examples below.

They used a Crank--Nicolson plus Strang splitting scheme with ghost-node Robin closure $u_{-1}=u_1+2\,dx\,(qu_0-h)$ and an absorbing layer. Frequencies were measured from the  trace by a zero-padded FFT with peak interpolation. They compared forced and unforced runs on the same grid, thus removing the $O(dx^2)$ shift of the computed frequencies.

\emph{One wall soliton}, $q=0.5$, $A=1$, $h=\varepsilon e^{-t/2}$ or resonant $h=i\varepsilon e^{-t/2}e^{it}$. The amplitude $A'$ was measured from the frequency, the trace amplitude and the mass; the three agree to $10^{-4}$.
A scan over $q\in\{-0.5,0,0.75,0.95\}$ with $\varepsilon=0.01$ gave measured-to-predicted ratios between $0.97$ and $1.01$.

\emph{Two wall solitons}, $q=0.3$, $\eta=(1.6,1.0)$, so $\omega=(2.56,1)$. They used $dx=0.02$, $dt=0.002$, $t\le450$, and measured frequencies on $[150,450]$.
Halving $\varepsilon$ roughly halves the deviation of the ratios from $1$, as expected for an $O(\varepsilon^2)$ remainder. A single forcing moves the two frequencies independently, and here in opposite directions. The late trace is again quasi-periodic, with both frequencies shifted.

\emph{Two wall solitons}, $q=-0.3$, $\eta =(0.8, 0.5)$, $h=\varepsilon(1-e^{-t})(1+t)^{-p}, \qquad t\le400 ,$ $\varepsilon=0.05$. Different values of $p$. Same phenomena: no new zeros, no non-imaginary zeros.
For larger $\varepsilon=0.3$ and $p=0.5$  more solitons can appear, but there are always at least two (sometimes more)  imaginary eigenvalues.

\section{Concluding remark}
The above elementary discussion is contrasted with the analysis of the Dirichlet problem. In that case  (see \cite{AK} and \cite{AK2}) one can show that some kind of integrability does persist when perturbing the linearisable
Dirichlet data (zero data). 

\section{Acknowledgements}
We are grateful for several discussions with Claude (Anthropic).

\end{document}